\documentclass[12pt]{article}
\usepackage{amsfonts,amsmath,amssymb}

\newtheorem{theorem}{Theorem}
\newtheorem{lemma}{Lemma}
\newtheorem{proposition}{Proposition}
\newtheorem{remark}{Remark}
\newtheorem{example}{Example}
\newtheorem{definition}{Definition}
\newtheorem{corollary}{Corollary}
\newlength{\rig}
\newlength{\rigg}
\newlength{\hei}
\newcommand{\dokaend}{\hfill$\square$ \vskip6truept}
\newcommand{\ind}{{\rm ind}}

\DeclareMathOperator{\im}{Im}
\DeclareMathOperator{\Vect}{Vect}
\DeclareMathOperator{\Ell}{Ell}
\DeclareMathOperator{\Mat}{Mat}
\DeclareMathOperator{\Con}{Con}

\title{Elliptic operators on manifolds \\
with singularities and
$K$-homology}
\author{Anton Savin}
\date{\today}

\begin{document}
\maketitle \abstract{ Elliptic operators on smooth compact manifolds are classified by
K-homology. We prove that a similar classification is valid also for manifolds with simplest
singularities: isolated conical points and edges. The main ingredients of the proof of these
results are: Atiyah-Singer difference construction in the noncommutative case and Poincare
isomorphism in K-theory for (our) singular manifolds. As an application we give a formula in
topological terms for the obstruction to Fredholm problems on manifolds with edges.\vspace{3mm}

}

%\tableofcontents

\section*{Introduction}
\addcontentsline{toc}{section}{Introduction}

It is well known \cite{Ati4,Kas3,BaDo1} that elliptic operators on a smooth closed manifold $M$
are classified by the $K$-homology group of the manifold. Namely, there is an isomorphism
\begin{equation}\label{atihomol}
\Ell(M)\simeq K_0(M),
\end{equation}
where $\Ell(M)$ is the group of stable homotopy classes of elliptic pseudodifferential operators
on the manifold, while  $K_0(M)$ is the (even) $K$-homology group, considered as the generalized
homology theory dual to Atiyah--Hirzebruch $K$-theory (see~\cite{White1}).

Both sides of  \eqref{atihomol} are well-defined even if the manifold or the operators have
singularities. More precisely, the right hand side is defined for finite complexes, while the
left hand side can be defined at least for simplest singularity types: isolated points, edges,
corners (see e.g. \cite{Schu1}, \cite{MaMe3}, and also \cite{Nis2}). It turns out that
\eqref{atihomol} holds for manifolds with singularities. As the main result of this paper, we
prove that the isomorphism is valid for manifolds with simplest singularities.

Let us  explain the scheme of obtaining the  isomorphism. We start from the
smooth case. Here the homotopy classification is given already by the Atiyah--Singer
difference construction
\begin{equation}\label{difc1}
\chi: \Ell(M) \stackrel\simeq\longrightarrow K_c(T^*M).
\end{equation}
Here $K_c$ denotes $K$-theory with compact support. Having in mind $K$-homological terminology,
we note that we can apply Poincar\'e isomorphism in $K$-theory on the cotangent bundle to the
latter group:
\begin{equation}\label{difc2}
p:K_c(T^*M)\stackrel\simeq\longrightarrow K_0(M)
\end{equation}
and obtain the desired isomorphism \eqref{atihomol}.

We now consider the case of singular $M$. Can one apply the same approach to
prove the desired isomorphism? It turns out that this scheme can be  applied
modulo certain modifications.

Namely, the symbols of operators on singular manifolds define
\emph{noncommutative algebras}. Hence, it is natural to work in the $K$-theory
of algebras.  Then the topological $K$-group $K_c(T^*M)=K^0_c(T^*M)$ of the
cotangent bundle in Eq.~\eqref{difc1} is replaced by the $K_0$-group of a
certain algebra (we denote it by $\mathcal{A}$) associated with the Calkin
algebra of our operators and the algebra that defines the bundles, in sections
of which our operators are acting. Hence, it remains to construct the
isomorphism
\begin{equation}\label{poi1}
p: K_0(\mathcal{A})\longrightarrow K_0(M),
\end{equation}
--- an analog of Poincar\'e isomorphism \eqref{difc2}, and we give here its construction on
the example of simplest singularities like the conical point and the edge.
Thus, we obtain the stable homotopy classification for elliptic operators on
manifolds with such singularities.

The homotopy classification \eqref{atihomol} has numerous applications. Let us
now note one  of them. In elliptic theory on manifolds with singularities, it
is well known \cite{Schu1}, \cite{NSScS4}, \cite{Nis1} that often an operator
elliptic in the interior part of the manifold may not be Fredholm. Moreover,
the Fredholm property can not be achieved even if we allow homotopies of the
operator near the singularity. (This phenomenon goes back to the theory of
boundary value problems, where there exists the \emph{Atiyah--Bott obstruction}
\cite{AtBo2} to elliptic boundary conditions for elliptic operators on
manifolds with boundary.) We will show that the homotopy classification enables
one to write explicit formula for the obstruction. This formula gives the same
result as in the previous computations in \cite{AtBo2}, \cite{SaSt11}.
\vspace{3mm}

Let us now briefly describe the contents of the paper.  In the first section we
give the abstract statement of the homotopy classification problem for elliptic
operators. Then we formulate the main results of the paper --- homotopy
classifications for simplest classes of manifolds with singularities: conical
points and edges. We would like to mention to simplify the presentation for
manifolds with edges we use the algebra of pseudodifferential operators with
discontinuous symbols introduced in \cite{SaSt11}. The proofs of the
classifications are given in the fourth section. They are based on the
generalization (in Section~3) of the Atiyah--Singer difference construction
\cite{AtSi1} to the noncommutative case and Mayer--Vietoris arguments. Section
5 contains some applications of the homotopy classification. We conclude the
paper with an appendix, where some basic facts about analytic $K$-homology are
summarized.

Let us mention some papers related to the results we obtain. The homotopy
classification for operators on manifolds with one  singular point was obtained
in  \cite{NScS5}. The Poincar\'e duality in terms of operator algebras on
manifolds  with isolated singularities was studied in a recent paper
\cite{DeLe1}. The  classification of some operators on smooth manifolds with
corners in terms of $K$-homology  is obtained in \cite{MePi2}. Added in proof.
The homotopy classification of edge-degenerate elliptic operators on manifolds
with edges was obtained in \cite{NSScS14} using the methods of the present
paper.

The author is grateful to Prof.~B.~Yu.~Sternin for a number of valuable remarks he made, while
the paper was in preparation.  The results of the paper were reported at the International
Conference ``Operator algebras on singular manifolds'' Potsdam, Germany, March 2003,  Workshop
"Index problems" Paris, April 2004 and other meetings. The work was partially supported by RFBR
grants NN~02-01-00118, 02-01-00928.

\section{Homotopy classification problem}\label{par1}

\textbf{Abstract operators in subspaces}. Consider a pair of algebras $ \mathcal{A}_0\subset
\mathcal{A} $ inside the algebra $\mathcal{ B}(H)$ of bounded operators on a separable Hilbert
space $H$. We will assume for simplicity that  $\mathcal{A}$ contains the ideal of compact
operators $\mathcal{K}=\mathcal{K}(H)$.

We call $\mathcal{A}$ the \emph{algebra of operators}, its Calkin algebra
$\mathcal{A}/\mathcal{K}$ ---  the \emph{algebra of symbols}, while algebra $\mathcal{A}_0$
--- \emph{algebra of functions}. By $\im A$ we denote the range of mapping $A$.

\begin{definition}
{\em {A linear operator}
\begin{equation}
\label{subsa1} D:\im P\longrightarrow \im Q
\end{equation}
is called \emph{operator in subspaces} (see \cite{SaSt1}) determined by the
pair
 $P,Q\in \Mat(\mathcal{A}_0,n)$ of matrix $n$ by $n$ projections ($P^2=P, Q^2=Q$) with entries in
 $\mathcal{A}_0$ and matrix operators $D\in \Mat(\mathcal{A},n)$
satisfying the inclusion $ D(\im P)\subset \im Q$ (algebraically this is equivalent to:
$QDP=DP$).

Two operators in subspaces are \emph{equal}, if they define the same restrictions
 \eqref{subsa1}. }
\end{definition}
This definition is an abstract analog of the notion of operator acting in sections of vector
bundles.
\begin{example}\label{odin1}
\emph{Consider the pair $\mathcal{A}_0=C^\infty(M),$ $\mathcal{A}={\Psi(M)}$,
where ${\Psi(M)}$ denotes the algebra of pseudodifferential operators of order
zero on a closed smooth manifold $M$. Both algebras are considered as
subalgebras of operators acting on  $L^2(M)$. In this case, projections with
entries in  $C^\infty(M)$ define sections of vector bundles over $M$,  while
the operator in subspaces \eqref{subsa1} is merely an operator acting in vector
bundle sections.}
\end{example}

\begin{definition}
\emph{ Operator $D:\im P\to \im Q$ is \emph{elliptic}, if there exists an almost-inverse
operator, i.e.  operator $D':\im Q\to \im P$, $D'\in\Mat(\mathcal{A},n)$ such that the
compositions $DD'$ and $D'D$ give (modulo compact operators) identity mappings
in subspaces  $\im Q$ and $\im P$.}
\end{definition}

\noindent \textbf{Stable homotopies. Group $\Ell(\mathcal{A}_0,\mathcal{A})$}. Two elliptic
operators $(D_0,P_0,Q_0)$  and $(D_1,P_1,Q_1)$ are  \emph{homotopic}, if they can be connected by
a continuous family of elliptic operators:\footnote{Here an operator family $D_t:\im P_t\to\im
Q_t$ (the source and the target spaces depend on $t$) is {\em continuous} if the projections
defining the spaces vary continuously and the family  $D_t$ is  a restriction of some continuous
family $\widetilde{D}_t\in\Mat(\mathcal{A},n)$.}
$$
D_t:\im P_t\to \im Q_t,\qquad t\in [0,1],\qquad  P_t,Q_t\in \Mat(\mathcal{A}_0,n),
$$
\begin{definition}
\emph{Two elliptic operators are \emph{stably homotopic}, if they become homotopic after addition
of some  trivial operators to each of them. Here \emph{trivial operators} are identity operators
of arbitrary matrix order (i.e., triples $(1,P,P)$ with an arbitrary projection $P$).}
\end{definition}

It is easy to check that stable homotopy is an equivalence relation. By
$\Ell(\mathcal{A}_0,\mathcal{A})$ we denote the set of elliptic operators,
modulo stable homotopies. This set is an Abelian group with respect to the
direct sum of operators in subspaces. The inverse element is given by the
almost-inverse operator and the zero is represented by trivial operators.

The main problem we consider in this paper is the problem of computing the group $\Ell$ for some particular
classes of algebras. We refer to this problem as the
\emph{homotopy classification problem for elliptic operators}  generated by the pair $({\cal A}_0,{\cal A})$
\vspace{2mm}

\noindent \textbf{Symbols.} Similar to operators one can introduce the notion of
 \emph{symbol in subspaces}
\begin{equation}\label{symbols}
(\sigma,P,Q).
\end{equation}
Here $P,Q$ are projections in matrix algebras over $\mathcal{A}_0$ as above,
while the matrix operator $\sigma\in \Mat(\mathcal{A}/\mathcal{K},n)$ has
components in the Calkin algebra. As usual it is easy to prove that the symbol
allows one to reconstruct the operator uniquely up to compact operators, and
the existence of a homotopy between two  operators is equivalent to the
existence of homotopies of their symbols. Therefore, classification of elliptic
operators up to stable homotopy is equivalent to a similar classification of
elliptic symbols. The latter is technically simpler.

\section{Main results}\label{par2}

\subsection{Operators on manifolds with conical points}

\textbf{Manifolds}. The simplest manifold with singularities is the \emph{
manifold $\mathcal{M}$ with a conical point}. It can be obtained if we identify
all the boundary points of a closed compact smooth manifold  $M$ with nonempty
boundary $\partial M$
$$
\mathcal{M}=M/\partial M.
$$
This space has only one point at which it is not smooth. This is the point, to which the
whole initial boundary is collapsed. Near the singularity the space looks like a neighborhood of
the vertex of the cone with base  $\partial M$
$$
K_{\partial M}=\partial M\times \overline{\mathbb{R}}_+/\partial M\times\{0\}.
$$

\noindent\textbf{Algebra of operators}. On $\mathcal{M}$, we consider the algebra
$\Psi(\mathcal{M})$ of cone pseudodifferential operators of order zero (e.g., see
\cite{NaSt12,EgSc1})
$$
D:L^2(\mathcal{M})\longrightarrow L^2(\mathcal{M}),
$$
acting in weighted Sobolev space of weight and order zero. These spaces will be denoted by
 $L^2(\mathcal{M})$.

Let us recall that the Calkin algebra\footnote{We shall assume that the
algebras of pseudodifferential operators contain the ideal of compact
operators.} $\Psi(\mathcal{M})/\mathcal{K}$ consists of compatible pairs
\begin{equation}\label{calkin1}
\Psi(\mathcal{M})/\mathcal{K}= \left\{ (\sigma,\sigma_c) \in C^\infty(S^*M)\oplus \Psi_p(\partial
M) \bigl | \;\; \sigma|_{\partial S^*M}=\text{ smbl } \sigma_c \right\}.
\end{equation}
\begin{itemize}
\item the first component $\sigma(D)$ is the  \emph{principal symbol}
--- function on the cosphere bundle $S^*M$ of the manifold, smooth up to the boundary
$\partial S^*M$;

\item the second component $\sigma_c(D)(p)$
--- the \emph{conormal symbol}, is a parameter-dependent family  with  $p\in\mathbb{R}$
(see \cite{AgVi1}) of zero order pseudodifferential operators on the boundary
$\partial M$. The algebra of parameter-dependent families is denoted by
$\Psi_p(\partial M)$.
\end{itemize}
The compatibility condition in \eqref{calkin1} involves the symbol mapping for
parameter-dependent families:
$$
{\rm smbl}:\Psi_p(\partial M)\longrightarrow C^\infty(S(T^*\partial M\times\mathbb{R})).
$$
Finally, we assume in \eqref{calkin1} that an identification $S^*M|_{\partial
M}=S(T^*\partial M\times\mathbb{R})$ is chosen.

We refer the reader to the cited papers for explicit formulas which allow one
to compute both symbols $\sigma(D),\sigma_c(D)$ for a given operator $D$ and,
conversely, to construct an operator starting from a compatible pair of symbols
$(\sigma,\sigma_c)$.

\vspace{3mm}

\noindent\textbf{Homotopy classification.} Denote by $\Ell(\mathcal{M})$ the group of stable
homotopy classes of elliptic operators associated with the inclusion
 $C^\infty(M)\subset\Psi(\mathcal{M})$:
$$
\Ell(\mathcal{M})=\Ell(C^\infty(M),\Psi(\mathcal{M})).
$$
In terms of elliptic operators one can say that this group is generated by elliptic operators
$$
D: L^2(\mathcal{M},E) \longrightarrow L^2(\mathcal{M},F),
$$
acting on sections of arbitrary vector bundles $E,F\in \Vect(M)$.

The  above description of the Calkin algebra $\Psi(\mathcal{M})/\mathcal{K}$ shows that the
commutator $[D,f]$  with function $f\in C^\infty(M)$ is compact if the function is constant on
the boundary. By continuity this gives the compactness of the commutator for an arbitrary
continuous function on the (compact) singular space $\mathcal{M}$
$$
[D,f]\in \mathcal{K},\qquad \text{ for all }  f\in C(\mathcal{M}).
$$
It follows (see Section~\ref{sechomol}) that an arbitrary elliptic operator $D$
defines an element in analytic $K$-homology of the singular space $\mathcal{M}$
(such elements were first defined in \cite{Has1}). We denote this element by
$$
\chi(D)\in K_0(\mathcal{M}).
$$

\begin{theorem}\label{thmain}
The mapping $D\mapsto \chi(D)$ defines an isomorphism of groups:
$$
\Ell(\mathcal{M}) \stackrel{\chi}\simeq K_0(\mathcal{M}).
$$
\end{theorem}

\noindent\textbf{Manifolds with several conical points.} Similarly, one can
consider manifold $\mathcal{M}$ with  $N>1$ distinct conical points. In other
words, when the singular space  is obtained from a manifold with boundary such
that the connected components of the boundary  are arranged into $N$
nonintersecting groups:
$$
\partial M=\Omega_1\sqcup \Omega_2\ldots \sqcup \Omega_N,
$$
and each group $\Omega_i$ is then collapsed to a point $pt_i$. Analytically, this means that we
consider operators with  $N$ conormal symbols. The isomorphism of Theorem \ref{thmain} remains
true in this case.

\begin{remark}
\emph{The homotopy classification in the case of one singular point was
obtained in \cite{NScS5}. It was shown that the group  $\Ell(\mathcal{M})$ is
isomorphic to the direct sum $K^0_c(T^*M)\oplus \mathbb{Z}$.
Theorem~\ref{thmain} gives the same answer in this case. Indeed:
$$
K_0(\mathcal{M})=K_0(M/\partial M)\simeq K_0(M,\partial M )\oplus \mathbb{Z}\simeq
K^0_c(T^*M)\oplus \mathbb{Z}.
$$
Here we use Poincar\'e isomorphism on manifolds with boundary \cite{Kas3}.}
\end{remark}

\subsection{Operators with discontinuous symbols in fiber bundles}\label{disc2}

In this subsection we study the classification problem for a class of elliptic
operators on fibered manifolds. These operators will be used in the next
subsection to obtain the homotopy classification on manifolds with fibered
boundary.

\noindent\textbf{Operators on a fibered manifold}. Consider a smooth locally trivial fiber bundle
$$
\pi:Y\to X.
$$
Assume that the base $X$ and the fiber (denoted by $\Omega$) are compact closed
manifolds.

Consider the algebra generated by pseudodifferential operators of order zero on $Y$ and families
of pseudodifferential operators of order zero acting in the fibers, see \cite{SaSt11}. Denote
this algebra by   $\Psi(Y,\pi)$. Its elements  can be considered also as operators  on  $Y$ with
symbols having discontinuities of a special form that we now describe. (Symbols discontinuous in
covariables  were studied earlier in \cite{Pla4}.)

Denote by $\overline{S^*Y\setminus \pi^*S^*X} $ the compactification of the complement
$S^*Y\setminus \pi^*S^*X$, for which a sequence $(x_i,\omega_i,\xi_i,\eta_i)$ converges to a
point in  $\pi^*S^*X$, if in addition to the usual convergence one also has the convergence for
the quotients $\eta_i/|\eta_i|$ (here $x,\omega$ are the coordinates along the base and the fiber, while
 $\xi,\eta$
--- are the dual coordinates). It is easy to show that $\overline{S^*Y\setminus
\pi^*S^*X} $ is a manifold with boundary.\footnote{The boundary fibers over $Y$, with fiber over
a point $(x,\omega)$ equal to the product of the cospheres $S^*_xX\times S^*_\omega\Omega$.}
This boundary is diffeomorphic to the complement $S^*Y\setminus U_{\pi^*S^*X}$ to a collar neighborhood
$U_{\pi^*S^*X}$ of the submanifold $\pi^*S^*X$.

Now the Calkin algebra of  $\Psi(Y,\pi)$ can be described (see~\cite{SaSt11}) as the subalgebra in
the direct sum of algebras:
\begin{equation}
\label{disc1} \Psi(Y,\pi)/\mathcal{K}\subset C^{\infty}(\overline{S^*Y\setminus \pi^*S^*X})\oplus
C^{\infty}(S^*X,\Psi(\Omega)),
\end{equation}
determined by the compatibility condition:
\begin{equation}\label{sogl}
\Psi(Y,\pi)/\mathcal{K}=\left\{(\sigma ,\widehat{\sigma})\quad\bigl|
\quad\sigma|_{\partial\overline{S^*Y\setminus \pi^*S^*X}}=\mathrm{smbl}\;
\widehat{\sigma}\right\}.
\end{equation}
For an operator $A\in \Psi(Y,\pi)$ the first component  $\sigma(A)$ is the
\emph{principal symbol} of the operator. The principal symbol is a function on
$S^*M$ smooth except $\pi^*S^*X$. The second component $\widehat{\sigma}(A)$
--- is called the \emph{operator symbol}. The operator symbol is a function on the cosphere
bundle of the base with values in  the algebras $\Psi(\Omega)$ of
pseudodifferential operators in the fibers.

Consider the pair of algebras $C^\infty(X,\Psi(\Omega))\subset \Psi(Y,\pi)$. The corresponding
group of stable homotopy classes of elliptic operators  is denoted by $\Ell(Y,\pi)$. It is
generated by operators
$$ D:\im P_1\longrightarrow \im P_2,
$$
acting on the fibered manifold, where the projections
$P_{1,2}:L^2(Y,E_{1,2})\to L^2(Y,E_{1,2})$ are induced by continuous families
of pseudodifferential operators in the fibers. We refer the reader to
\cite{SaSt11} for a detailed study of these operators, in particular, for the
proof of \eqref{sogl}.\vspace{2mm}

\noindent{\bf Homotopy classification.} The properties of the Calkin algebra
 $\Psi(Y,\pi)/\mathcal{K}$ show that for an elliptic operator
$D\in\Ell(Y,\pi)$ the commutator $[D,f]$  is compact provided $f$ is constant in the fibers of
 $\pi$. Therefore, operator $D$ defines an element
$$
\chi(D)\in K_0(X)
$$
in analytic  $K$-homology of the base.
\begin{theorem}\label{thfib}
One has an isomorphism:
\begin{equation}
\label{fibered1} \Ell(Y,\pi)\stackrel\chi\simeq K_0(X).
\end{equation}
\end{theorem}

\begin{remark}\label{zama1}
\emph{ A similar result holds in the case, when the base of the fiber bundle has a nonempty
boundary. In this case we consider the algebra generated by the operators which near the boundary
act as operators of multiplication by functions $f\in C^\infty(Y)$. The corresponding group
$\Ell(Y,\pi)$ will be isomorphic in this case to the group $K^0_c(T^*\overset \circ X)\simeq
K_0(X)$, where $\overset \circ X$ denotes the interior part.  The same classification can be
obtained if one considers more general operators, which near the boundary are defined by
operator-valued functions  $f\in C^\infty(X,\Psi(\Omega))$. This result can be obtained by the
same method.}
\end{remark}

\subsection{Operators on manifolds with fibered boundary}

\noindent\textbf{Manifolds with fibered boundary.} Let  $M$ be a compact manifold with boundary,
and the boundary is the total space of a locally--trivial fiber bundle with fiber $\Omega$. We
will fix some extension of the fibration to some collar neighborhood of the boundary $U_{\partial
M}\simeq
\partial M\times [0,1)$ and denote the extended projection by
$$
\pi:\partial M\times [0,1]\to X\times [0,1].
$$

Let  $\mathcal{M}$ be the manifold with singularities, obtained from $M$ by the identification of
points in the fibers of
 $\pi$. A homotopy equivalent space can be obtained if we identify only the points in the fibers
 at the boundary $\partial M$. This space is referred to as a
\emph{manifold with edge} $X$.

The aim of the present subsection is to describe an elliptic theory, whose
homotopy classification produces the group $K_0(\mathcal{M})$.\vspace{2mm}

\noindent\textbf{Operator algebra.} Denote by $\Psi(M,\pi)\subset
\mathcal{B}(L^2(M))$ the algebra generated by:
\begin{itemize}
\item the usual pseudodifferential operators of order zero over the interior part of the
manifold;

\item in a neighborhood  $\partial M\times [0,1)$ of the boundary by operators from the algebra
$\Psi(\partial M\times [0,1),\pi)$ (see Remark~\ref{zama1} of the previous section);

\item compact operators.
\end{itemize}
In other words, the operators from the algebra $\Psi(M,\pi)$ are the usual pseudodifferential
operators far from the boundary, and in a neighborhood of the boundary are operators
corresponding to the fibration $\pi$.\vspace{2mm}

\noindent\textbf{Homotopy classification.} Let $\Ell(M,\pi)$ be the group of stable homotopy
classes of elliptic operators corresponding to the embedding of algebras $C^\infty(M)\subset
\Psi(M,\pi)$.

The symbolic calculus of the usual pseudodifferential operators and operators in the fibration
$\pi$ show that for an  operator $D\in\Psi(M,\pi)$, and $f\in C^\infty(M)$ the commutator $[D,f]$
is compact, if the function is constant in the fibers of  $\pi$, i.e., it is induced by a
continuous function on the manifold with edge. Hence (see Section~\ref{sechomol}), an elliptic
operator $D$ from $\Ell(M,\pi)$ defines an element in $K$-homology of the singular space
$\mathcal{M}$. Denote this element as
$$
\chi(D)\in K_0(\mathcal{M}).
$$
\begin{theorem}\label{thfb1}
One has
\begin{equation}
\Ell(M,\pi) \stackrel\chi\simeq   K_0(\mathcal{M}).
\end{equation}
\end{theorem}

\section{$\Ell$-theory and $K$-theory}

\subsection{Difference construction}

The aim of this subsection is to prove that the group  $\Ell(\mathcal{A}_0,\mathcal{A})$ is
isomorphic to the $K$-group of some algebra constructed from the pair
$(\mathcal{A}_0,\mathcal{A})$. The isomorphism generalizes the Atiyah--Singer difference
construction \cite{AtSi1}.

In addition to the conditions of Section~\ref{par1}, we will assume that
$\mathcal{A}$ and $\mathcal{A}_0$ are $C^*$-algebras. For definiteness we will assume that
 $\mathcal{A}_0$ is unital.

By $\Con_f$ denote the \emph{cone of a monomorphism} $f:A_0\to A$. We will also write
$\Con(A_0,A)$, when $f$ is the inclusion of a subalgebra. Recall that the cone is defined as the
subalgebra of the direct sum
 $A_0\oplus C_0([0,1),A)$:
\begin{equation}
\label{konus1} \Con_f=\left\{(a_0,a(t))\in {A}_0\oplus C_0\Bigl([0,1),A\Bigr) \;\Bigl|\;
f(a_0)=a(0)\right\}.
\end{equation}

\noindent\textbf{Difference element of an elliptic symbol.} Let us prove that
the symbol of an elliptic operator in subspaces
$$
D:\im P\to \im Q
$$
defines an element in the $K$-group $K_0(\Con(\mathcal{A}_0,\mathcal{A}/\mathcal{K}))$
of the cone of the natural mapping $\mathcal{A}_0\to \mathcal{A}/\mathcal{K}$. This element is
constructed as follows.

To shorten the notation, we denote the algebra $\Con(\mathcal{A}_0,\mathcal{A}/\mathcal{K})$ by
$\Con$, while the algebra with adjoined unit by $\Con^+$. The latter algebra also consists of
pairs $(a,a(t))$. The only difference is that for $\Con^+$ the function $a(t)$ at  $t=1$ can be
an arbitrary scalar:
$$
a(1)=\lambda Id,\qquad \lambda\in \mathbb{C}.
$$

Let us define a matrix projection $\mathcal{P}$ with entries in $\Con^+$:
\begin{equation}\label{proj1}
\mathcal{P}=(P,P_t)\in\Mat(\Con^+,N).
\end{equation}
The second term in this pair
$P_t\in\Mat(\mathcal{A}/\mathcal{K},N)$ is a family of projections defined for $t\in[0,1]$ as
\begin{equation}\label{povorot}
P_{t}=P\cos^2\varphi+Q\sin^2\varphi+\cos\varphi\sin\varphi(\sigma P+\sigma^{-1}Q), \quad
\varphi=\frac\pi 2t
\end{equation}
(here $\sigma\in A$ denotes the symbol of $D$, $\sigma^{-1}$ is the symbol of the almost-inverse
operator). The homotopy of projections \eqref{povorot} is called \emph{rotation homotopy}.\footnote{The
term {\em rotation} is motivated by the fact that in the case of orthogonal projections
 $P,Q$ and an isometric isomorphism $\sigma$ the vectors in the range of $P_t$
are obtained as the rotation  by angle $\pi t/2$ of  $v\in \im P$ towards $\sigma v\in \im Q$.}
The homotopy is well defined (i.e., for all  $t$ the operator $P_t$ is a projection), provided
the ranges of $P$ and $Q$ are orthogonal. The orthogonality is valid, if, for instance,
the projections have the following block-diagonal form
$$
P=\left(%
\begin{array}{cc}
  \boxed{*} & 0 \\
  0 & 0 \\
\end{array}%
\right),\quad Q=\left(%
\begin{array}{cc}
  0 & 0 \\
  0 & Id_n \\
\end{array}%
\right),\qquad P,Q\in\Mat(\mathcal{A}_0,m+n)
$$
(an arbitrary pair $(P,Q)$ can be reduced to this form by stable homotopies).

Let $\mathcal{P}_0\in\Mat(\Con^+,m+n)$ be the projection with components $(Q,Q)$ (cf.~\eqref{proj1}). The element
\begin{equation}\label{qdiff1}
\chi(\sigma(D))=[\mathcal{P}]-[\mathcal{P}_0]\in K_0(\Con),
\end{equation}
is called the  \emph{difference element} of the elliptic operator $D$. \vspace{2mm}

\noindent \textbf{Difference construction}
\begin{theorem}\label{abstr}
The mapping $D\mapsto \chi(\sigma(D))$ induces a  group isomorphism \emph{(}difference
construction\emph{)}
\begin{equation}\label{diffk}
\Ell(\mathcal{A}_0,\mathcal{A})\stackrel{\chi}\simeq
K_0(\Con(\mathcal{A}_0,\mathcal{A}/\mathcal{K}) ).
\end{equation}
\end{theorem}

\begin{example}
\emph{ On a smooth manifold, the Calkin algebra
$\overline{\Psi(M)}/\mathcal{K}$ (here $\overline{\Psi(M)}$ is the norm closure
in $\mathcal{B}(L^2(M))$) is isomorphic to the algebra $C(S^*M)$ of continuous
functions on the cosphere bundle of the manifold, the cone of the embedding
$C(M)\subset C(S^*M)$ is isomorphic to the algebra $C_0(T^*M)$ of functions on
the cotangent bundle vanishing at infinity. In this case the theorem gives the
Atiyah--Singer isomorphism: the group of stable homotopy classes of elliptic
pseudodifferential operators on
 $M$ is isomorphic to the  $K$-group  $K_c^0(T^*M)$. }
\end{example}

\noindent\emph{Proof.} Let us prove that the difference construction is well
defined. Actually, this follows from the definition of the groups $\Ell$ and
$K_0$. Indeed, the difference element of a trivial operator is zero in the
$K$-group, while a homotopy of elliptic operators induces a homotopy of the
corresponding projections $\mathcal{P}$ and $\mathcal{P}_0$ in \eqref{qdiff1}.

The invertibility of the difference homomorphism is proved by first constructing the left inverse
mapping to $\chi$ and then proving the surjectivity of $\chi$. These two facts prove the desired
isomorphism. Let us prove these facts.

1. Let us construct the inverse mapping
\begin{equation}\label{inva3}
\chi^{-1}:K_0(\Con)\to \Ell(\mathcal{A}_0,\mathcal{A}).
\end{equation}
Consider an element $[\mathcal{P}]-[\mathcal{Q}]\in K_0(\Con)$, where $\mathcal{P},\mathcal{Q}$
are projections in matrix algebras over $\Con^+$. We can assume that $\mathcal{Q}=Id_n\oplus
0_m$. Then for the projection $\mathcal{P}=(P,P_t)$ we define
$$
\chi^{-1}([\mathcal{P}]-[\mathcal{Q}])\stackrel{\rm def}=\left[U:\im P\longrightarrow \im
P_1\right],
$$
where the symbol of an elliptic operator $U$ is defined as the value at $t=1$ of the solution of the
Cauchy problem:
\begin{equation}\label{defo1}
\left\{%
\begin{array}{c}
  \dot{u}_t=[\dot{P}_t,P_t]u_t, \\
  u_{0}=1. \\
\end{array}%
\right.
\end{equation}
This procedure requires the smoothness in $t$ of the homotopy. Nonetheless, this is sufficient to
define the mapping \eqref{inva3}. This follows
from the fact that the subalgebra in  $\Con$, corresponding to smooth families of projections
$P_t$, is local  (in the sense of  \cite{Bla1}) and, thus, defines the same $K$-group as the algebra
$\Con$.

2. The equality $\chi^{-1}\circ\chi=Id$ is obtained by a direct computation, since
\eqref{defo1} can be solved explicitly for the  projection defined by
\eqref{povorot}. The details of the computation are left to the reader.

3. To prove the surjectivity of $\chi$ one has to show that an arbitrary element
\begin{equation}
\label{mudiff1} [\mathcal{P}]-[\mathcal{Q}]\in K_0(\Con)
\end{equation}
can be represented in the form
 \eqref{qdiff1}, for  a rotation homotopy \eqref{povorot} defined by
some elliptic symbol $\sigma$.

Let us deform the pair $\mathcal{P},\mathcal{Q}$ to the desired form. To this end, it is enough
to assume that the second projection is trivial: $\mathcal{Q}=(Q,Q)$, as in \eqref{qdiff1}. In
this case $Q$ is determined by $\mathcal{P}=(P,P_t)$, or, more precisely, by the value of its
second component at $t=1$. Hence, it will be enough to construct only homotopies of
$\mathcal{P}$, without mentioning explicitly the homotopy of $\mathcal{Q}$.

We construct the desired homotopy $\mathcal{P}=(P,P_t)$ in two steps.

First,  we deform the homotopy  to obtain a  homotopy
 $P_t$, which for $t\in[1/2,1]$ is the rotation homotopy from
 $P_{1/2}$ towards $P_1$, while the projections  $P_t$ for all $t\le
1/2$ are orthogonal to  $P_1$. This can be accomplished by a stabilization (passing to matrices
twice bigger as the original ones) and by superposition of the original homotopy with the
rotation homotopy, which connects the projections $P_1\oplus 0$ and $0\oplus P_1$.

Second, we extend the rotation homotopy obtained previously from the half-interval
 $[1/2,1]$ to the whole
$[0,1]$. To this end we deform the homotopy of projections   $P_t$ with deformation parameter
$\varepsilon\in[0,1/2]$ according to the formula
$$
P_{\varepsilon,t}=\left\{
\begin{array}{cc}
P_t, & \text{ for  }  t\le\varepsilon, \vspace{2mm}\\
\begin{array}{c}
\text{ rotation by angle }\frac{\displaystyle (t-\varepsilon)\pi}{\displaystyle (1-\varepsilon)2}
\vspace{1mm}\\
\text{ from } P_\varepsilon \text{ towards } P_1,
\end{array}
& \text{ for  }   t\ge\varepsilon,
\end{array}
\right.
$$
where in the second case the rotation between the projections with  a given homotopy between them is
defined by substituting in \eqref{povorot} the solution of the Cauchy problem \eqref{defo1}. By
construction, for $\varepsilon=0$ we obtain the rotation homotopy between  $P_0$ and $P_1$.
Therefore, the initial element  \eqref{mudiff1} is indeed a difference element
for some elliptic operator.

This proves that the difference construction is surjective. \dokaend

\begin{remark}\label{famil1}
\emph{(families of operators) One can generalize the difference construction to the case of
continuous \emph{families of elliptic operators}. This generalization is done by standard
techniques. Since this will be used later (see Remark~\ref{remka}), we will formulate the result:
families of elliptic operators corresponding to a  pair $(\mathcal{A}_0,\mathcal{A})$ and
parametrized by a compact space  $X$ are classified by the  $K$-group of the algebra $C(X,\Con)$
of continuous $\Con$-valued functions on $X$.}
\end{remark}

\begin{remark}\label{remka}
\emph{Operator families can be used to define the \emph{odd} elliptic group corresponding to the
group  $K_1$ of the mapping cone. Namely, let  $\Ell_1(\mathcal{A}_0,\mathcal{A})$ be the group
of stable homotopy classes of elliptic families  $\{D_\varphi\}_{\varphi\in\mathbb{S}^1}$ such that
 for $\varphi=0$ $D_\varphi$ is trivial (i.e., the operator at $\varphi=0$ has components in the
algebra  $\mathcal{A}_0$). A family is trivial, if $D_\varphi$ is trivial for all  $\varphi$. It
is easy to express the odd groups in terms of the even ones:
$$
\Ell_1(\mathcal{A}_0,\mathcal{A})\simeq \Ell(\Sigma \mathcal{A}_0,\Sigma\mathcal{A}).
$$
(Here $\Sigma A$ is the suspension $C_0((0,1),A)$ of $A$.) Hence, applying Theorem ~\ref{abstr}
we get the desired isomorphism
$$
\Ell_1(\mathcal{A}_0,\mathcal{A})\simeq K_1(\Con(\mathcal{A}_0,\mathcal{A}/\mathcal{K})).
$$
Here we use the suspension isomorphism.}
\end{remark}

\subsection{Mapping cone $K$-theory exact sequence}
One has an exact sequence
\begin{equation}\label{as1}
\ldots\to K_1(\mathcal{A}/\mathcal{K})\to K_0(\Con)\to K_0(\mathcal{A}_0)\longrightarrow
K_0(\mathcal{A}/\mathcal{K})\to K_1(\Con)\ldots
\end{equation}
induced by the short exact sequence of algebras
\begin{equation}
\label{mapa1} 0\to \Sigma (\mathcal{A}/\mathcal{K}) \longrightarrow \Con \longrightarrow
\mathcal{A}_0 \to 0.
\end{equation}
The elements in  \eqref{as1}  admit natural descriptions in terms of elliptic operators:

\begin{itemize}
\item the group  $K_1(\mathcal{A}/\mathcal{K})$ corresponds to the subclass of matrix elliptic
operators (i.e., the projections  $P,Q$ are identity maps) and is formed by the stable homotopy
classes of such operators;

\item the mapping $K_1(\mathcal{A}/\mathcal{K})\to K_0(\Con)$ in \eqref{as1}
is induced by  the embedding of matrix operators in the class of  operators in arbitrary
subspaces (i.e., subspaces defined by possibly  nontrivial projections);

\item the mapping $K_0(\Con)\to K_0(\mathcal{A}_0)$ takes an operator in subspaces to the
difference of its projections;

\item the boundary mappings $K_*(\mathcal{A}_0)\longrightarrow K_*(\mathcal{A}/\mathcal{K})$
are induced by the algebra homomorphism
$f:\mathcal{A}_0\to\mathcal{A}/\mathcal{K}$ (e.g., see \cite{CuSk1});

\item the mapping $K_0(\mathcal{A}/\mathcal{K})\to K_1(\Con)$ takes a projection
$P\in\Mat(\mathcal{A}/\mathcal{K},n)$ to the family $Pz+Id_n-P$,  where $z=e^{i\varphi}$ is the
coordinate on the circle.
\end{itemize}

\section{Proofs of the homotopy classifications}

In this section we give the proofs of the theorems stated in Section~\ref{par2}.

\subsection{Isolated singularities}

In this subsection we give the proof of Theorem~\ref{thmain} on the homotopy classification of
elliptic operators on a manifold with a conical point.

1. Instead of proving the classification of elliptic operators, we shall classify elliptic
symbols (according to Section~\ref{par1} the two problems are equivalent).

2. Furthermore, the classification corresponding to the pair of algebras
 $C^\infty(M)\subset \Psi(\mathcal{M})$ coincides with the classification
for the closure of this pair with respect to the operator norm:  $C(M)\subset
\overline{\Psi(\mathcal{M})}$. Let us note that similarly to the original algebra, the Calkin
algebra of the closure consists of compatible pairs
$$
\overline{\Psi(\mathcal{M})}/\mathcal{K}= \left\{ (\sigma,\sigma_c)\in C(S^*M)\oplus
\overline{\Psi_p(\partial M)} \Bigl | \;\;\sigma|_{\partial S^*M}=\text{ smbl } \sigma_c
\right\}.
$$
Here the closure $\overline{\Psi_p(\Omega)}$  of the algebra of parameter-dependent families is
considered with respect to the norm
$$
\sup_{p\in\mathbb{R}}\|\sigma_c(p)\|_{L^2(\partial M)\to L^2(\partial M)}.
$$
This description follows from the well-known estimates of the norm for operators
 $D\in\Psi(\mathcal{M})$ modulo sums with compact operators:
$$
\inf_{K\in \mathcal{K}}\|D+K\|_{L^2(\mathcal{M})\to L^2(\mathcal{M})}=
\mathrm{max}\bigl(\max_{(x,\xi)\in S^*M}|\sigma(D)(x,\xi)|,
\sup_{p\in\mathbb{R}}\|\sigma_c(D)(p)\| \bigr).
$$

3. We can apply the difference construction of the previous section (Theorem~\ref{abstr}) to the
embedding of  $C^*$-algebras $f:C(M)\to \overline{\Psi(\mathcal{M})}/\mathcal{K}$:
$$
\Ell(\mathcal{M})\simeq K_0(\Con_f).
$$
Thus, we will prove the isomorphism of the latter group and the $K$-homology group of the
singular manifold.

4. Let us embed $\chi:K_0(\Con_f)\to K_0(\mathcal{M})$ in the diagram:
\begin{equation}
\begin{array}{ccccccccc}
\label{di1}
\to K^1_c(T^*M)&\stackrel\partial\to &\mathbb{Z}&\to &K_0(\Con_f)&\to& K_c(T^*M)&\to& 0 \\
\downarrow & & \parallel & & \chi\downarrow & & \downarrow \\
\to K_1(M,\partial M)&\stackrel{\partial ''}\to& K_0(pt)&\to& K_0(\mathcal{M})&\to&
K_0(M,\partial M)&\to& 0.
\end{array}
\end{equation}
The bottom row of the diagram is the $K$-homology exact sequence for the pair $pt\subset
\mathcal{M}$,  while the upper row is the $K$-theory exact sequence for the ideal
$$
I\subset \Con_f
$$
of elements with zero principal symbol. For this ideal it is easy to obtain an isomorphism
 $I\simeq C_0((0,1)\times \mathbb{R},\mathcal{K})$, and for the quotient we have
$$
\Con_f/I\simeq C_0(B^*M\setminus S^*M)\simeq C_0(T^*M).
$$
(Here $B^*M$ is the unit ball bundle in the cotangent bundle.)

The vertical mappings in the diagram are induced by quantizations. More precisely, the elements
of the $K$-groups in the upper row are considered as difference elements for some elliptic
symbols, and then one associates operators to the symbols. In more detail, $\chi$ corresponds to
cone degenerate operators. The mappings $K^*_c(T^*M)\to K_*(M,\partial M)$ on the sides of the
diagram are well known (e.g., see \cite{Kas3}, \cite{BDT1}) and correspond to operators elliptic
on the interior of the manifold (let us recall that in the case of odd  $K$-homology groups one
considers only self-adjoint elliptic operators). Such operators almost-commute with all functions
vanishing on the boundary. All the vertical mappings, except possibly $\chi$, are isomorphisms
(e.g., see \cite{Kas3}).

5. Let us compute the boundary mapping $\partial$. To this end we interpret the group
 $K^1_c(T^*M)\simeq K_c(T^*M\times\mathbb{R})$ as the group generated by symbols
 $\sigma(x,\xi,\lambda)$ elliptic with parameter $\lambda\in
\mathbb{R}$. Here the coordinates $x$ in a neighborhood of the boundary are separated into the
tangential and normal  $x=(y,t)$, the dual variables are denoted by $\xi=(\eta,p)$. Then the
restriction of $\sigma(x,\xi,\lambda)$ to the boundary can be considered as the family of symbols
on the boundary elliptic with parameters $(p,\lambda)\in \mathbb{R}^2$. Therefore, the index of this
family is a number
$$
\ind\; \sigma\left(y,0,-i\frac\partial{\partial y},p,\lambda\right) \in K_c(\mathbb{R}^{2})\simeq
\mathbb{Z}.
$$
This is well defined, since parameter dependent ellipticity implies that the operator family
$\sigma\left(y,0,-i\partial/{\partial y},p,\lambda\right)$ is Fredholm everywhere and invertible
for large values of the parameters  (e.g., see \cite{Shu1}).

\begin{lemma} The boundary $\partial[\sigma]$ of an  element  $[\sigma]\in K^1_c(T^*M)$  is
expressed by the formula
\begin{equation}\label{bou1}
\partial[\sigma]=\ind \;\sigma\left(y,0,-i\frac\partial{\partial y},p,\lambda\right)\in
K_c(\mathbb{R}^{2})\simeq \mathbb{Z}.
\end{equation}
\end{lemma}

\noindent\emph{Proof.} Consider a commutative diagram
\begin{equation}\label{coma1}
\begin{array}{ccccccccc}
 0 & \to & I & \to & \Con_f & \to & C_0(T^*M) & \to & 0 \\
 & & \parallel & & \downarrow & &  \downarrow i \\
 0 & \to & I & \to & \partial\Con_f  & \to & C_0(T^*M|_{\partial M}) & \to &
 0,
\end{array}
\end{equation}
here $\partial \Con_f$ is the algebra
$$
 C_0([0,1),\overline{{\Psi}_p(\partial M)})\cap \{\sigma_c(0) - \text{\rm
operator of multiplication}\}=\Con(C(\partial M),\overline{\Psi_p(\partial M)}).
$$
The diagram \eqref{coma1} and the naturality of the boundary mapping imply that the desired
boundary mapping  $\partial$ (it corresponds to the top row in \eqref{coma1}) is a composition
$$
\partial=\partial'i_*
$$
of the restriction  $i_*:K^*_c(T^*M)\to K^*_c(T^*M|_{\partial M})$ and the boundary mapping
$$\partial':K^*_c(T^*M|_{\partial M}) \to K_{*+1}(I)$$ corresponding to the bottom row.

We are interested in the boundary mapping $\partial'$ defined on the odd $K$-group. Hence, the
desired equality  \eqref{bou1} follows from the fact that the boundary mapping from the odd group
is actually the index  (see \cite{Bla1}). \dokaend

6. We can apply the 5-lemma to the diagram and thus prove the desired isomorphism
 $K_0(\Con_f)\simeq
K_0(\mathcal{M})$, once we know the diagram is commutative. Therefore, the proof will be
completed if we obtain the following lemma.
\begin{lemma}
The diagram \emph{(\ref{di1})} commutes.
\end{lemma}
\emph{ Proof. } A. Let us verify the commutativity of the square
\begin{equation*}
\begin{array}{ccc}
K^1_c(T^*M)&\to &\mathbb{Z}\\
\downarrow & & \parallel \\
K_1(M,\partial M)&\stackrel{\partial ''}\to& K_0(pt).
\end{array}
\end{equation*}
The boundary mapping $\partial''$ can be represented as the composition of the boundary mapping
$K_1(M,\partial M)\to K_0(\partial M)$ and the index mapping $K_0(\partial M)\to
K_0(pt)=\mathbb{Z}$, i.e.  $\partial''$ coincides with the mapping in the upper row of the
diagram \eqref{di1} (see also \cite{MePi1}).

B. Commutativity of the square in the center of \eqref{di1}. The range of  $\mathbb{Z}\to
K_0(\Con_f)$ is defined by elliptic operators of the form $1+\mathbf{G}$, where $\mathbf{G}$ is
an operator with principal symbol equal to zero. We have to prove that such operators give
Fredholm modules equivalent to a Fredholm module with a very simple module structure: the action
of a function $f$ is determined by the product with the value $f(pt)$ at the cone vertex. This is
accomplished in two steps:
\begin{enumerate}
\item (removal of the smooth part) without changing the element in $K$-homology, we can
restrict the operator $1+\mathbf{G}$ to a neighborhood of the conical point  (the
original operator and its restriction are stably equivalent);

\item (homotopy of the module structure) in a neighborhood of the conical point the module
structure  $f(x),u(x)\mapsto f(x)u(x)$ is reduced by a homotopy to $f(x),u(x)\mapsto f(pt)u(x)$.
The homotopy corresponds to the simplest rescaling: $f(x\varepsilon)u(x)$.
\end{enumerate}

C. Commutativity of the right most square follows immediately from the definition of $K$-homology
for noncompact spaces (i.e.,   $K$-theory for nonunital algebras), see  \cite{BDT1}.\dokaend

Thus, Theorem~\ref{thmain} is proved for the case, when there is only one singular point. For
several conical points the ideal $I$ of operators with zero principal symbol
is a direct sum of ideals $C_0(\mathbb{R}^2,\mathcal{K})$ at each  singular point,
and we obtain
$$
K_0(I)\simeq \mathbb{Z}^N.
$$
The scheme of the proof remains the same in this case and we shall not repeat it here.

\subsection{Fibered manifolds}

In this subsection we prove Theorem~\ref{thfib} on the homotopy classification of elliptic
operators on a fiber bundle.

1. As in the previous subsection we consider the classification for elliptic symbols.
Furthermore, we classify continuous symbols.  The inclusion of norm-closures of our algebras is
denoted by
$$
f:C(X,\overline{\Psi(\Omega)})\longrightarrow \overline{\Psi(Y,\pi)}/\mathcal{K}.
$$
The Calkin algebra of the closure $\overline{\Psi(Y,\pi)}$ is the subalgebra of compatible pairs
of continuous symbols as in  \eqref{sogl}:
$$
\overline{\Psi(Y,\pi)}/\mathcal{K}\subset C(\overline{S^*Y\setminus \pi^*S^*X})\oplus
C(S^*X,\overline{\Psi(\Omega)}).
$$
This follows from the equality of the norms of operators modulo compact terms, and the maximums
of the absolute values of symbols (see \cite{SaSt11})
\begin{equation}
\inf_{K\in \mathcal{K}}\|A+K\|_{L^2(Y)\to L^2(Y)}=\max\left( \max_{S^*Y}|\sigma(A)|,
\max_{S^*X}\|\widehat{\sigma}(A)\|_{L^2(\Omega)\to L^2(\Omega)} \right).
\end{equation}

2. We embed $\chi$ in a commutative diagram:
\begin{equation}
\label{treug1}
\begin{array}{ccc}
K^0_c(T^*X)\!\! &\!\!\longrightarrow\!\!  &\!\! K_0(\Con_f) \vspace{2mm}\\
& \simeq\searrow\quad& \downarrow \chi \vspace{2mm}\\
& & K_0(X).
\end{array}
\end{equation}
Here the horizontal mapping is induced by the embedding of the operators on the base $X$ in the
set of operators on the fibration. This embedding is defined by realizing the usual finite-rank
vector bundles on $X$ as the ranges of families of finite rank projections in the $L^2$-spaces on
the fibers.

3. We will obtain  the desired isomorphism property for $\chi$, if we prove that the horizontal
mapping of the triangle is an isomorphism. This is proved in the following proposition.
\begin{proposition}
The inclusion of the ideal $I_X\subset \Con_f$ of elements with zero principal symbol
induces the isomorphism of $K$-groups
$$
K_0(\Con_f)\simeq K_0({I}_X)\simeq K^0_c(T^*X).
$$
\end{proposition}
\emph{Proof}.  The quotient $\Con_f/I_X\simeq C_0(\overline{T^*M\setminus \pi^* T^*X})$
corresponds to the product $[0,\infty)\times \partial \left(\overline{T^*Y\setminus \pi^*
T^*X}\right)$, i.e., is contractible and has trivial $K$-groups. Thus, the short exact sequence
of the pair  $I_X\subset \Con_f$ induces the desired isomorphism of the  $K$-groups of the algebra and
its ideal.

The isomorphism  $K_0({I}_X)\simeq K^0_c(T^*X)$ follows from the fact that the ideal $I_X$
consists of functions on  $T^*X$ with values in compact operators in the fibers. \dokaend

\subsection{Manifolds with fibered boundary}

In this subsection we prove Theorem~\ref{thfb1} about the homotopy classification of
elliptic operators on manifolds with fibered boundary.

1. The difference construction gives an isomorphism $\Ell(M,\pi)\simeq K_0(\Con)$, where $\Con$
denotes the cone of  the embedding:
$$
C(M)\longrightarrow \overline{\Psi(M,\pi)}/\mathcal{K}.
$$
Here $C(M)$ and $ \overline{\Psi(M,\pi)}$ are the norm-closures of the smooth algebra of functions and the algebra
of operators.

2. The groups $\Ell$ and $K$- will not change, if we replace the latter map by the monomorphism
$$
f:\widetilde{C}_0(M)\longrightarrow \overline{\Psi_0(M,\pi)}/\mathcal{K}.
$$
Here $\Psi_0(M,\pi)\subset \Psi(M,\pi)$ is the subalgebra of operators with symbols equal to zero
on $\partial M$, and the algebra of ``functions'' $\widetilde{C}_0(M)$ is defined as:
$$
\widetilde{C}_0(M)=\left\{ (u,v) \in C(M\setminus U_{\partial M})\oplus C(X\times
(0,1],\overline{\Psi(\Omega)})
 \bigl|
u_{t=1}=v_{t=1} \right\}.
$$
Here  $t\in[0,1]$ is a normal coordinate in a neighborhood $U_{\partial M}$ of the boundary,
$\{t=0\}$ is the equation of the boundary. This choice of the algebras corresponds to the class
of operators introduced in Remark~\ref{zama1}.

3. We will prove the isomorphism for the $K$-groups $K_0(\Con_f)$ and $K_0(\mathcal{M})$.

4. The groups are embedded in the diagram:
\begin{equation}\label{atom}
\begin{array}{ccccccccc}
K^1_c(T^*M)&\to &K_0(\mathcal{A}_0) & \to & K_0(\Con_f) & \to &  K_c^0(T^*M) &
\to & K_1(\mathcal{A}_0) \\
\chi_M \downarrow\simeq& & \chi_X\downarrow\simeq & & \chi_{\mathcal{M}}\downarrow & &
\chi_M\downarrow\simeq & & \chi_X\downarrow\simeq\\
K_1(M\setminus \partial M)\!&\!\! \stackrel\partial\to \!\!\!&\!\! K_0(X\times [0,1])\! &\! \to
\!&\! K_0(\mathcal{M}) \!&\! \!\to\!\! & K_0({M}\setminus \partial M) \!&\!\!
\stackrel{\partial}\to\!\! & K_1(X\times[0,1]).
\end{array}
\end{equation}
Let us define this diagram.
The bottom row here is the $K$-homology exact sequence of the pair
 $X\times[0,1]\subset \mathcal{M}$. The top row is induced by the exact sequence of algebras
\begin{equation}\label{seqex}
 0\longrightarrow \mathcal{A}_0 \longrightarrow \Con_f\stackrel{j}\longrightarrow
C_0(T^*M_0)\longrightarrow 0,
\end{equation}
here $j$ is the restriction mapping to the closed subset  $M_0=M\setminus U_{\partial M}$, where
we deal with the usual scalar symbols. The kernel of  $j$ is an ideal in $\Con_f$ denoted by
$\mathcal{A}_0$. One can see that this ideal is the mapping cone of the homomorphism:
$$
f_{(0,1)}:C_0(\partial M\times (0,1))\to \overline{\Psi(\partial M\times (0,1),\pi)}/\mathcal{K}.
$$
In particular, according to Theorem~\ref{abstr} its $K$-group $K_0(\mathcal{A}_0)\simeq K_0(X)$
classifies  elliptic operators corresponding to the fibration  $\pi:\partial M\times[0,1]\to
X\times [0,1]$ (see Subsection~\ref{disc2}).

The vertical mappings in \eqref{atom} --- quantizations $\chi$ with lower indices --- are defined
in terms of elliptic operators:  $\chi_M$ corresponds to operators on the interior of $ M$,
$\chi_X$ ---  to operators on the fibration $\pi$.

All the vertical mappings, except $\chi_{\mathcal{M}}$, are isomorphisms. Therefore, to
prove that
 $\chi_{\mathcal{M}}$ is an isomorphism using the 5-lemma, it is enough to prove the commutativity of the diagram.

5. Let us show that the left- and rightmost squares in \eqref{atom} commute.  The commutative
diagram
\begin{equation}\label{maps1}
\begin{array}{ccccc}
\partial M\times [0,1] & \subset & M &  \longrightarrow & M/\partial M \\
\pi\downarrow\quad & & \downarrow & & \parallel\\
X \times[0,1]& \subset & \mathcal{M} & \to & M/\partial M
\end{array}
\end{equation}
implies that the boundary mapping $\partial$ in \eqref{atom} is a composition
\begin{equation}\label{alo1}
K_*(M,\partial M) \longrightarrow K_{*+1}(\partial M)\stackrel{\pi_*}\longrightarrow K_{*+1}(X)
\end{equation}
of the restriction to the boundary and the change of the module structure (direct image mapping
$\pi_*$).  Applying Poincar\'e isomorphism (see \cite{BaDo1}, \cite{Kas3}), we can pass to the
topological $K$-groups, and rewrite  \eqref{alo1} as
$$
K^*_c(T^*M)\longrightarrow K^{*}_c(T^*M|_{\partial M})=K^{*+1}_c(T^*\partial M)
\stackrel{\pi_!}\longrightarrow K^{*+1}_c(T^*X),
$$
here $\pi_!$ is the direct image mapping in  $K$-theory induced by $\pi$.

It remains to compare this mapping with the boundary mapping $ K^*_c(T^*M)\longrightarrow
K_{*+1}(\mathcal{A}_0)$ in the top row of \eqref{atom}. The latter mapping is also a
composition
$$
K^*_c(T^*M)\stackrel{i_*}\longrightarrow K_c^*(T^*M|_{\partial
M})\stackrel{\partial'}\longrightarrow K_{*+1}(\mathcal{A}_0)
$$
of the restriction to the boundary
 $\partial M$ and the boundary mapping in $K$-theory (denoted by  $\partial
'$) corresponding to the bottom row of the following diagram
\begin{equation}
\label{zeta}
\begin{array}{ccccccccc} 0 \longrightarrow &\Con_{f_{(0,1)}} &\longrightarrow & \Con_f &
\longrightarrow & C_0(T^*M)& \longrightarrow & 0\\
 & \parallel & &  \downarrow & & \downarrow i \\
0 \longrightarrow &\Con_{f_{(0,1)}}
 &\longrightarrow & \Con_{f_{(0,1]}} & \longrightarrow &
C_0(T^*M|_{\partial M})& \longrightarrow & 0.
\end{array}
\end{equation}
Here $f_{(0,1]}$ is the restriction of $f$ to $(0,1]\times \partial M$.

Let us compute $\partial'$. To this end, we rewrite the sequence in terms of symbols on the
boundary. By $\partial f:C(X,\overline{\Psi(\Omega)}) \to \overline{\Psi(\partial
M,\pi)}/\mathcal{K}$ denote the inclusion of algebras associated with the fibered boundary.

The bottom row sequence in \eqref{zeta} can be represented in the form \eqref{mapa1}. Namely, it
is obtained
$$
0\to \mathcal{A}_0=\Sigma^{2}\Con_{\partial f}\to \Sigma\Con_{h}\to \Sigma C_0(T^*\partial M)\to
0
$$
as the suspension\footnote{The suspension common to all terms of the sequence uses the
conormal variable.} of the sequence of the form \eqref{mapa1} induced by the homomorphism
$$
h:C_0(T^*\partial M)\to \Con_{\partial f}.
$$
The latter homomorphism is induced by the embedding of the continuous symbols in the class of
discontinuous symbols.

Thus, according to \eqref{as1} the boundary $\partial'$  in $K$-theory
$$
K_*(C_0(T^*\partial M))\stackrel {\partial '}\longrightarrow K_*(\Con_{\partial f})
$$
coincides with the induced mapping $h_*$.

Therefore, the commutativity of the extreme squares in \eqref{atom} reduces to the commutativity
of the following diagram
$$
\begin{array}{ccc}
K_*(C_0(T^*\partial M)) & \stackrel {h_*}\longrightarrow & K_*(\Con_{\partial f}) \vspace{1mm}\\
\chi_{\partial M} \downarrow \quad & & \quad\downarrow  \chi_X \vspace{1mm}\\
K_*(\partial M) & \stackrel {\pi_*}\longrightarrow & K_*(X).
\end{array}
$$
It is commutative, since  $\partial M$ is closed and smooth.

6. The commutativity of the two squares in the middle of
 \eqref{atom} is proved easily.

Applying the  5-lemma to \eqref{atom} we end the proof of Theorem \ref{thfb1}.

\section{Applications and remarks}

\noindent \textbf{Obstructions to Fredholm problems for elliptic operators}.
Let $\mathcal{M}$ be a manifold with edge $X$. Choose a diffeomorphism
$\mathcal{M}\setminus X\simeq \overset \circ M$ of the smooth part and the
interior of a smooth compact  manifold $M$ with boundary. Let $\pi:\partial
M\to X$ be the natural projection.

Let $D$ be an operator on $\mathcal{M}$ with elliptic interior principal
symbol. As we know, this condition does not guarantee the Fredholm property of
$D$. Then the following question naturally arises: can one make $D$ Fredholm by
some modification of lower order terms of the operator?

It appears that the homotopy classification enables one to answer a slightly
weaker version of this question, when instead of modifying lower order terms
one can make stable homotopies of the principal symbol of the operator.
\begin{corollary}
Let $[D]\in K_0(\mathcal{M}\setminus X)$ be an element defined by an operator
$D$ with elliptic interior principal symbol. Then this element can be lifted to
an element in $K_0(\mathcal{M})\simeq \Ell(\mathcal{M})$ defined by some
operator elliptic on the entire space $\mathcal{M}$ if and only if
$$
\partial [D]=0,
$$
where $\partial:K_0(\mathcal{M}\setminus X)\longrightarrow K_1(X)$ is the
boundary map of the exact sequence of the pair $X\subset\mathcal{M}$.
\end{corollary}
This readily follows from the exact sequence $\to K_0(\mathcal{M})\to
K_0(\mathcal{M}\setminus X)\stackrel\partial \to K_1(X)\to$.

In specific situations this formula gives many obstructions known previously.

1. For $X=\partial M$, i.e. when $\mathcal{M}=M$ is a manifold with boundary,
the vanishing of $\partial [D]\in K_1(\partial M)$  is equivalent to the
vanishing of the Atiyah--Bott \cite{AtBo2} element $[\sigma(D)|_{\partial
T^*M}]\in K^0(T^*M|_{\partial M})\simeq K^1(T^*\partial M)$. The equivalence is
given by Poincare duality $K_*(\partial M)\simeq K^*(T^*\partial M)$).

2. If $\mathcal{M}$ is a general manifold with edges, then  the application of
Poincare duality to $\partial[D]$ gives the element
$$
\pi_![\sigma(D)|_{\partial T^*M}]\in K^1(T^*X).
$$
This formula for the obstruction was obtained in \cite{SaSt1} using a different
method.

3. For isolated singularities ($\dim X=0$) the obstruction can be nontrivial
only in the self-adjoint setting, i.e., for $[D]\in K_1(\mathcal{M}\setminus
X)$. In this case the obstruction reduces to the so called deficiency indices,
e.g. see~\cite{Les1}, which are the obstructions to self-adjoint extensions of
symmetric operators.

Let us finally note that the equality $\partial[D]=0$ can be effectively
verified modulo torsion.
\begin{proposition}
$\partial[D]=0\in K_1(X)\otimes \mathbb{Q}$ if and only if
\begin{equation}
\label{pari1} \langle [\sigma(D)|_{\partial T^*M}],\pi ^*a\rangle=0
\end{equation}
for all elements $a\in K^1(X)$. Here
$$
\langle,\rangle:K^1(T^*\partial M)\times K^1(\partial M)\to K^0(T^*\partial
M)\stackrel\ind\longrightarrow \mathbb{Z}
$$
is the natural index pairing.
\end{proposition}
\emph{Proof}. By the naturality of the boundary mapping and the index pairing
we get:
\begin{equation}\label{notin1}
\langle \partial[D], a\rangle=\langle \pi_*\partial'[D], a\rangle=\langle
\partial'[D], \pi^*a\rangle,
\end{equation}
where $\partial':K_0(\mathcal{M}\setminus X)=K_0(M\setminus \partial M)\to
K_1(\partial M)$ is the boundary map for the pair $\partial M\subset M$.

By Poincare duality \eqref{notin1} is equal to $\langle [\sigma(D)|_{\partial
T^*M}], \pi^*a\rangle$. The equality $\langle \partial[D], a\rangle=\langle
[\sigma(D)|_{\partial T^*M}], \pi^*a\rangle$ proves the desired statement,
since the pairings are nondegenerate (Poincare duality). \dokaend

\noindent \textbf{$K$-groups of $C^*$-algebras of pseudodifferential
operators}. The problem of computing the $K$-groups of norm-closures of
algebras of pseudodifferential operators on manifolds with singularities was
considered in \cite{MeNi2}. Let us show that the answer to this problem can be
given using the homotopy classification.

First, we claim that the embedding $f:C(M)\to \overline{\Psi(\mathcal{M})}/\mathcal{K}$ induces a
monomorphism in $K$-theory if $\mathcal{M}$ has no closed smooth components. This follows from
the existence of a nonsingular vector field on $M$: the vector field defines a left inverse
$\overline{\Psi(\mathcal{M})}/\mathcal{K}\to C(S^*M)\to C(M)$ to $f$, and, hence, the $6$-term
exact sequence splits:
$$
K_i(\overline{\Psi(\mathcal{M})}/\mathcal{K})\simeq K^{i}(M)\oplus K_{i+1}(\Con_f).
$$
Second,  the exact sequence of the pair $\mathcal{K}\subset\overline{\Psi(\mathcal{M})}$ gives
$$
K_0(\overline{\Psi(\mathcal{M})})\simeq K_0(\overline{\Psi(\mathcal{M})}/\mathcal{K})\qquad
K_1(\overline{\Psi(\mathcal{M})})\simeq
\ker\left(K_1(\overline{\Psi(\mathcal{M})}/\mathcal{K})\stackrel{\ind}\to \mathbb{Z}\right).
$$
Combining these two results and the homotopy classification
$K_{*}(\Con_f)\simeq K_*(\mathcal{M})$, we obtain
\begin{corollary}
$$
K_i(\overline{\Psi(\mathcal{M})})\simeq K^i(M)\oplus \widetilde{K}_{i+1}(\mathcal{M}),
$$
where $\widetilde{K}_*$ is the reduced $K$-homology group.
\end{corollary}

\noindent \textbf{Classification, Poincare duality in $K$-theory. Relation to
index theory}. In terms of $K$-theory for $C^*$-algebras the homotopy
classification gives an isomorphism
$$
K_*(\mathcal{A})\stackrel\simeq\to K^*(C(\mathcal{M}))
$$
of $C^*$-algebraic $K$-groups of opposite variance. In $K$-theory such
isomorphisms usually appear in pair with isomorphisms
$$
K^*(\mathcal{A})\stackrel\simeq\to K_*(C(\mathcal{M}))
$$
of dual theories and one refers to $\mathcal{A}$ and $C(\mathcal{M})$ as
\emph{Poincare dual} algebras \cite{Con1}, \cite{Kas3}.

Suppose for the moment that such a dual pair of isomorphisms exists. Denote by
$x\in K^0(\mathcal{A})$ the preimage of the identity in
$K_0(C(\mathcal{M}))=K^0(\mathcal{M})$. Then for an elliptic operator $D$ on
$\mathcal{M}$ we obtain:
$$
\ind D=\langle[D],1\rangle=\langle[\sigma(D)],x\rangle,\qquad \text{where }
[\sigma(D)]\in K_0(\mathcal{A}),
$$
by Poincare duality. The expression of the index in this form is important,
because it contains only the principal symbol of $D$. If element $x$ is
explicitly given then one can try to pass by Chern--Connes character to cyclic
cohomology \cite{Con1} and compute the corresponding class ${\rm ch} x\in
HC^*(\mathcal{A})$ to obtain an index formula.

It would be very interesting to find a realizations of the element $x$ by some
operator of geometric origin (some important steps in this direction can be
found in \cite{DeLe1} and \cite{NRSS1}). Let us finally note that for a smooth
manifold the element $x$ is given by the Dirac operator on $T^*M$. Thus, the
problem is to realize the element $x$ by a Dirac like operator.

Let us note finally that the  Poincare duality on manifolds with singularities
is useful in the theory of index defects, see \cite{SaSt10}, where it is used
to obtain index defect formulas. These topics for manifolds with singularities
will be discussed elsewhere.

\section{Appendix. Analytic $K$-homology}\label{sechomol}

In this Appendix we recall for the readers convenience the  basic definitions
from the theory of analytic $K$-homology (see \cite{HiRo1}, \cite{Bla1},
\cite{SoTr1}), which are used throughout the paper. Let $X$ be a compact
Hausdorf topological space.\vspace{2mm}

\noindent{\bf Even Fredholm modules and the group $K_0(X)$.}
\begin{definition}
\label{d8} {\em An \emph{even Fredholm module} over space $X$ is a pair $(F,H)$, consisting of a
bounded operator
$$
F:H\to H,
$$
acting on the $\mathbb{Z}_2$-graded separable Hilbert space
$$
H=H_0\oplus H_1,
$$
with components $H_{0,1}$ equipped with the structure of $*$-modules over the $C^*$-algebra
$C(X)$.  Denote the homomorphism defining the module structure by $\phi: C(X)\to \mathcal{B}(H)$.
We assume that the operator is odd relative to the grading, and for an arbitrary function
$f\in C(X)$ we have
\begin{equation}
\label{kompa} \phi(f)(F-F^*)\in {\cal K}(H), \quad \phi(f)(F^2-1)\in {\cal K}(H),\quad[F,\phi(f)]\in {\cal K}(H),
\end{equation}
where ${\cal K}(H)$ is the ideal of compact operators in $H$.}
\end{definition}

The  {\em even $K$-homology group} (with respect to direct sum) denoted by $K_0(X)$ is obtained
from Fredholm modules, if we identify the modules by the equivalence relation
--- \emph{stable homotopy}.\footnote{In \cite{Bla1}, chapter VIII,
one can find a number of other equivalence relations for Fredholm modules, which for a large class
of $C^*$-algebras are equivalent to the relation we use in this paper.}Let us describe this equivalence
relation.

Two Fredholm modules are \emph{isomorphic}, if the corresponding $C(X)$-modules $H$ are
isomorphic, and the operators $F$ transform one into another under the isomorphism of spaces. Two
modules $(F_1,H)$ and $(F_0,H)$ are \emph{homotopic}, if they can be connected by a family of
modules $(F_t,H)$ such that the operator  families  $F_t$ and $\phi_t(f)$ (for fixed $f\in C(X)$)
are strongly continuous. A module is \emph{trivial} if it has each  expression in \eqref{kompa}
equal to zero (for any $f\in C(X)$). Finally, two modules are \emph{stably homotopic}, if their
direct sums with some trivial modules  are homotopic.\vspace{2mm}

\begin{remark}
\emph{The definition of the $K$-homology groups still makes sense if one replaces the commutative
algebra $C(X)$  by much more general $C^*$-algebras. However, in this paper we use only the
algebra $C(X)$ of continuous functions on a compact set or the algebra $C_0(M\setminus \partial
M)$ of functions vanishing on the boundary of a smooth compact manifold $M$. The $K$-homology
group of the latter nonunital algebra is called the {\em relative $K$-homology group} and denoted
by $K_0(M,\partial M )$.}
\end{remark}

\noindent{\bf Odd Fredholm modules and the group $K_1$.}
\begin{definition}
An odd Fredholm module {\em over the space  $X$ is a pair $(F,H)$. Here $F$ is a bounded
operator
$$
F:H\to H,
$$
acting on a separable Hilbert space $H$, which is a $*$-module over the $C^*$-algebra $C(X)$ and for an
arbitrary function
 $f\in C(X)$ one has
\begin{equation}
\label{kompa1} \phi(f)(F-F^*)\in {\cal K}(H), \quad \phi(f)(F^2-1)\in {\cal K}(H),\quad[F,\phi(f)]\in {\cal K}(H).
\end{equation}
Here $\phi$ is the module structure.
}
\end{definition}
Thus, in contrast with the even case, the grading is not required. The set of stable homotopy
classes of odd Fredholm modules over  $X$ is denoted by $K_1(X)$.\vspace{2mm}

\noindent\textbf{Elliptic operators define elements of the groups $K_*$}. The usual elliptic
(pseudo)differential operators on a smooth closed manifold define elements in $K$-homology.
Namely, if $D$ is an elliptic pseudodifferential operator in sections of some bundles $E$ and $F$
over a smooth closed manifold  $M$, then it defines a bounded operator on $L^2$-spaces
\begin{equation}\label{normali1}
D'=(P_{\ker D}+D^*D)^{-1/2}D:L^2\left( M,E\right) \longrightarrow L^2\left( M,F\right).
\end{equation}
(Here $P_{\ker D}$ is the orthogonal projection on the kernel of $D$.) Both $L^2$-spaces of
sections are modules over the algebra   $C\left( M\right) $. To construct a Fredholm module
starting from $D'$ one considers the matrix operator
\[
F=\left(
\begin{array}{cc}
0 & D'^{*} \\
D'& 0
\end{array}
\right)
\]
as a self-adjoint operator in the naturally  $\Bbb{Z}_2$-graded $C\left(
M\right) $-module $H=L^2\left( M,E\right)\oplus L^2\left( M,F\right) $.
Operator $F$ is odd relative to the grading. The corresponding element in
$K$-homology is denoted by
\begin{equation}
\label{khom} \chi(D)\stackrel{def}=[F,H]\in K_0(M).
\end{equation}

Odd Fredholm modules can be obtained from \emph{self-adjoint} elliptic operators. Namely, if
 $A$
is an elliptic self-adjoint operator  in sections of a vector bundle
 $E$ on a closed manifold, then the odd Fredholm module is defined using the operator
\begin{equation}
\label{normali2} A'=(P_{\ker A}+A^2)^{-1/2}A:L^2\left( M,E\right)
\longrightarrow L^2\left( M,E\right) .
\end{equation}
The corresponding element in $K$-homology is denoted by
$$
\chi(A)\stackrel{def}=[A',H]\in K_1(M).
$$

%\bibliographystyle{unsrt}
%\bibliography{elliptic,algebras}

\begin{thebibliography}{10}

\bibitem{Ati4}
M.~F. Atiyah.
\newblock Global theory of elliptic operators.
\newblock In {\em Proc. of the Int. Symposium on Functional Analysis}, pages
  21--30, Tokyo, 1969. University of Tokyo Press.

\bibitem{Kas3}
G.~Kasparov.
\newblock Equivariant ${KK}$-theory and the {N}ovikov conjecture.
\newblock {\em Inv. Math.}, 91(1):147--201, 1988.

\bibitem{BaDo1}
P.~Baum and R.~G. Douglas.
\newblock ${K}$-homology and index theory.
\newblock In R.~Kadison, editor, {\em Operator Algebras and Applications},
  number~38 in Proc. Symp. Pure Math, pages 117--173. American Mathematical
  Society, 1982.

\bibitem{White1}
George~W. Whitehead.
\newblock Generalized homology theories.
\newblock {\em Trans. Amer. Math. Soc.}, 102:227--283, 1962.

\bibitem{Schu1}
B.-W. Schulze.
\newblock {\em Pseudodifferential Operators on Manifolds with Singularities}.
\newblock North--Holland, Amsterdam, 1991.

\bibitem{MaMe3}
R.~R. Mazzeo and R.~B. Melrose.
\newblock {P}seudodifferential operators on manifolds with fibred boundaries.
\newblock {\em Asian J. Math.}, 4:833--866, 1998.

\bibitem{Nis2}
V.~Nistor.
\newblock Singular integral operators on non-compact manifolds and analysis on
  polyhedral domains.
\newblock arXiv: math.AP/0402322, 2004.

\bibitem{NSScS4}
V.~Nazaikinskii, A.~Savin, B.-W. Schulze, and B.~Sternin.
\newblock { On the existence of elliptic problems on manifolds with edges}.
\newblock \emph{Russ. Acad. Sci. Doklady}, 69(2):231--234, 2004.

\bibitem{Nis1}
V.~Nistor.
\newblock An index theorem for gauge-invariant families: {T}he case of solvable
  groups.
\newblock {\em Acta Math. Hungarica}, 99(2):155--183, 2003.

\bibitem{AtBo2}
M.~F. Atiyah and R.~Bott.
\newblock The index problem for manifolds with boundary.
\newblock In {\em Bombay Colloquium on Differential Analysis}, pages 175--186,
  Oxford, 1964. Oxford University Press.

\bibitem{SaSt11}
A.~Savin and B.~Sternin.
\newblock {\em Boundary Value Problems on Manifolds with Fibered Boundary}.
\newblock Chalmers Institute of Technology, G{\"o}teborg, December 2001.
\newblock Preprint 75, available at arxiv.org/math.OA/0207179.

\bibitem{AtSi1}
M.~F. Atiyah and I.~M. Singer.
\newblock The index of elliptic operators {I}.
\newblock {\em Ann. of Math.}, 87:484--530, 1968.

\bibitem{NScS5}
V.~Nazaikinskii, B.-W. Schulze, and B.~Sternin.
\newblock {\em On the Homotopy Classification of Elliptic Operators on
  Manifolds with Singularities}.
\newblock Univ. Potsdam, Institut f{\"u}r Mathematik, Potsdam, Oktober 1999.
\newblock Preprint N 99/21.

\bibitem{DeLe1}
C.~Debord and J.-M. Lescure.
\newblock ${K}$-duality for pseudomanifolds with isolated singularities.
\newblock preprint math.OA/0212120, 2002.

\bibitem{MePi2}
R.~Melrose and P.~Piazza.
\newblock Analytic {$K$}-theory on manifolds with corners.
\newblock {\em Adv. in Math.}, 92(1):1--26, 1992.

\bibitem{NSScS14}
V.~Nazaikinskii, A.~Savin, B.-W. Schulze, and B.~Sternin.
\newblock {\em On the {H}omotopy {C}lassification of {E}lliptic {O}perators on
  {M}anifolds with {E}dges}.
\newblock Univ. Potsdam, Institut f{\"u}r Mathematik, Potsdam, 2004.
\newblock Preprint N 2004/16.

\bibitem{SaSt1}
A.~Yu. Savin and B.~Yu. Sternin.
\newblock Elliptic operators in even subspaces.
\newblock {\em Matem. sbornik}, 190(8):125--160, 1999.
\newblock English transl.: Sbornik: Mathematics {\bf 190}, N 8 (1999), p.
  1195--1228; arXiv: math/9907027.

\bibitem{NaSt12}
V.~E. Nazaikinskii and B.~Yu. Sternin.
\newblock Operator algebras on manifolds with isolated singularities.
\newblock {\em Differential Equations}, 39(1):92--104, 2003.

\bibitem{EgSc1}
Yu. Egorov and B.-W. Schulze.
\newblock {\em Pseudo-Differential Operators, Singularities, Applications}.
\newblock Birkh{\"a}user, Boston, Basel, Berlin, 1997.

\bibitem{AgVi1}
M.~Agranovich and M.~Vishik.
\newblock Elliptic problems with parameter and parabolic problems of general
  type.
\newblock {\em Uspekhi Mat. Nauk}, 19(3):53--161, 1964.
\newblock English transl.: Russ. Math. Surv. {\bf 19} (1964), N 3, p. 53--157.

\bibitem{Has1}
P.~Haskell.
\newblock Index theory of geometric {F}redholm operators on varieties with
  isolated singularities.
\newblock {\em $K$-Theory}, 1(5):457--466, 1987.

\bibitem{Pla4}
B.~A. Plamenevskii.
\newblock On pseudodifferential operators with discontinuities in the symbols
  with respect to momenta and coordinates.
\newblock {\em Dokl. Akad. Nauk}, 356(5):599--601, 1997.

\bibitem{Bla1}
B.~Blackadar.
\newblock {\em $K$-Theory for Operator Algebras}.
\newblock Number~5 in Mathematical Sciences Research Institute Publications.
  Cambridge University Press, 1998.
\newblock Second edition.

\bibitem{CuSk1}
J.~Cuntz and G.~Skandalis.
\newblock Mapping cones and exact sequences in ${KK}$-theory.
\newblock {\em J. Operator Theory}, 15(1):163--180, 1986.

\bibitem{BDT1}
P.~Baum, R.~G. Douglas, and M.~E. Taylor.
\newblock Cycles and relative cycles in analytic ${K}$-homology.
\newblock {\em J. Differ. Geom.}, 30(3):761--804, 1989.

\bibitem{Shu1}
M.~A. Shubin.
\newblock {\em Pseudodifferential Operators and Spectral Theory}.
\newblock Springer--Verlag, Berlin--Heidelberg, 1985.

\bibitem{MePi1}
R.~Melrose and P.~Piazza.
\newblock Families of {D}irac operators, boundaries and the $b$-calculus.
\newblock {\em J. of Diff. Geometry}, 46(1):99--180, 1997.

\bibitem{Les1}
M.~Lesch.
\newblock {\em Differential {O}perators of {F}uchs {Type}, {Conical}
  {S}ingularities, and {A}symptotic {M}ethods}, volume 136 of {\em
  Teubner--Texte zur Mathematik}.
\newblock B. G. Teubner Verlag, Stuttgart--Leipzig, 1997.

\bibitem{MeNi2}
R.~Melrose and V.~Nistor.
\newblock {$K$}-theory of {${\bf C}^*$}-algebras of $b$-pseudodifferential
  operators.
\newblock {\em Geom. Funct. Anal.}, 8(1):88--122, 1998.

\bibitem{Con1}
A.~Connes.
\newblock {\em Noncommutative geometry}.
\newblock Academic Press Inc., San Diego, CA, 1994.

\bibitem{NRSS1}
V.~Nazaikinskii, G.~Rozenblioum, A.~Savin, and B.~Sternin.
\newblock {G}uillemin {T}ransform and {T}oeplitz {R}epresentations for
  {O}perators on {S}ingular {M}anifolds.
\newblock In G.~Grubb B.~Booss-Bavnbek and K.~P. Wojciechowski, editors, {\em
  {S}pectral {G}eometry of {M}anifolds with {B}oundary}, Contemp. Math., 2005.
\newblock (in print).

\bibitem{SaSt10}
A.~Savin and B.~Sternin.
\newblock Index defects in the theory of nonlocal boundary value problems and
  the $\eta$-invariant.
\newblock {\em Sbornik:Mathematics}, 195(9), 2004.
\newblock Preliminary version at arXiv: math/0108107.

\bibitem{HiRo1}
N.~Higson and J.~Roe.
\newblock {\em Analytic {$K$}-homology}.
\newblock Oxford University Press, Oxford, 2000.

\bibitem{SoTr1}
Yu.~P. Solovyov and E.~V. Troitsky.
\newblock {\em {$C\sp *$}-algebras and elliptic operators in differential
  topology}, volume 192 of {\em Translations of Mathematical Monographs}.
\newblock American Mathematical Society, Providence, RI, 2001.

\end{thebibliography}
%\end{document}

%\vspace{1cm}

%\hfill {\em Moscow, Potsdam}

\end{document}